# From Finite Cayley Graphs to Growth of Infinite Groups

Tal Weissblat

Email: tal.weisblat@gmail.com

## Abstract

Graph neural networks (GNNs) have recently been shown to learn algebraic properties of finite groups from their Cayley graphs [1,2]. In this work, we investigate whether such models generalize to infinite finitely generated groups. Motivated by Gromov's theorem [3], a GNN is trained and validated exclusively on finite complete and truncated Cayley graphs, and then evaluated, without retraining, on truncated Cayley graphs of unseen infinite groups. The evaluation includes free abelian groups of various ranks, the discrete Heisenberg group, the infinite dihedral group, free groups, and direct products with both infinite abelian and finite groups. The results show strong generalization across these families, suggesting that finite Cayley graphs encode sufficient local geometric information to transfer to the infinite setting. Overall, this provides evidence that GNNs trained solely on finite groups can capture geometric features related to the growth of infinite finitely generated groups.

## 1. Introduction

Recent advances have demonstrated that graph neural networks (GNNs) can successfully learn algebraic properties directly from Cayley graphs of finite groups. In our previous work, we showed that GNNs can predict algebraic properties such as solvability and introduced a general framework for learning algebraic properties from finite Cayley graphs [1,2]. These results suggest that GNNs are capable of extracting meaningful algebraic and geometric information from graph representations of groups.

A fundamental question is whether the representations learned from finite groups can generalize beyond the finite setting. This question is closely related to one of the central themes of geometric group theory: the relationship between the geometry of Cayley graphs and the algebraic structure of infinite finitely generated groups. In particular, Gromov's celebrated theorem [3] states that an infinite finitely generated group has polynomial growth if and only if it is virtually nilpotent, establishing a profound connection between the large-scale geometry of its Cayley graph and its algebraic structure. This naturally raises the question of whether GNNs trained only on finite Cayley graphs can infer geometric information relevant to the growth of previously unseen infinite groups.

In this paper, we investigate whether GNNs trained exclusively on finite Cayley graphs can generalize to previously unseen infinite finitely generated groups. A GNN is trained and validated using only finite groups and is subsequently evaluated, without retraining, on truncated Cayley balls of infinite groups. The evaluation includes free abelian groups of several ranks, the discrete Heisenberg group, the infinite dihedral group, free groups, and direct products with both infinite abelian factors and finite cyclic and symmetric groups.

The experimental results demonstrate that the learned representations transfer successfully across this diverse collection of infinite groups. In particular, the trained models accurately distinguish groups of polynomial growth from groups of exponential growth while remaining robust under direct products with finite groups. These findings suggest that, when finite training groups provide a good approximation of the local geometry of infinite groups at the chosen truncation radius, GNNs are able to learn geometric features that generalize beyond the finite examples used during training.

To the best of our knowledge, this is the first experimental study investigating whether GNNs trained solely on finite Cayley graphs can generalize to infinite finitely generated groups for the purpose of predicting their growth behavior. The results provide evidence that graph representation learning can capture geometric information relevant to group growth, suggesting a promising connection between machine learning and geometric group theory.

The remainder of this paper is organized as follows. Section 2 reviews the necessary background on group growth, Cayley graphs, and the geometric aspects of finitely generated groups. Section 3 describes the experimental methodology, including the datasets, graph representations, and GNN architecture. Section 4 presents the experimental results, and Section 5 concludes with a discussion of the findings and their implications for the generalization of GNNs from finite to infinite groups.

## 2. Group Theory Background

This section briefly reviews the mathematical concepts used throughout the paper, including Cayley graphs, truncated Cayley balls, and the growth of finitely generated groups [4].

**Cayley graphs**. Let $G$ be a finitely generated group with a finite generating set $S$, where $S = S^{-1}$. The Cayley graph of $G$ with respect to $S$, denoted by $Cay(G, S)$, is the undirected graph whose vertices are the elements of $G$. Two vertices $g, h \in G$ are adjacent if and only if $h = gs$ for some generator $s \in S$. Equivalently, for every generator $s$, the graph contains an undirected edge between $g$ and $gs$. Cayley graphs provide a geometric representation of finitely generated groups by encoding the action of the generators as graph edges. Consequently, many algebraic and geometric properties of a group are reflected in the local and global structure of its Cayley graph.

**Truncated Cayley balls**. Let $G$ be a finitely generated group with generating set $S$, and let $R \geq 0$. The truncated Cayley ball of radius $R$ is the induced subgraph of $Cay(G, S)$ consisting of all vertices whose word distance from the identity is at most $R$, together with all edges between such vertices. Thus, it captures the local geometry of the Cayley graph within distance $R$ of the identity.

Throughout this work, finite groups are represented either by their full Cayley graphs or by truncated Cayley balls, depending on the experimental setting. Since the Cayley graph of an infinite finitely generated group is itself infinite, infinite groups are necessarily represented by truncated Cayley balls. Consequently, both finite and infinite groups are presented to the GNN as finite graphs, allowing the same GNN to be trained on finite groups and then evaluated directly on previously unseen infinite groups without retraining.

**Growth**. The growth of an infinite finitely generated group describes how the number of group elements within word distance $r$ of the identity increases as $r$ grows. Equivalently, it measures the growth of the balls in the Cayley graph of the group. Two of the most important growth types are polynomial growth and exponential growth [4].

A fundamental theorem due to Gromov states that an infinite finitely generated group has polynomial growth if and only if it is virtually nilpotent [3]. Consequently, determining the growth type of a group is closely related to understanding its underlying algebraic structure. This connection provides the mathematical motivation for the present work, which investigates whether GNNs can infer growth behavior from finite portions of Cayley graphs.

## 3. Methodology

This section describes the experimental methodology used to investigate whether GNNs trained exclusively on finite groups can generalize to unseen infinite finitely generated groups. The following figure provides an overview of the proposed framework. The remainder of this section describes the construction of the training, validation, and testing datasets, followed by the graph representation, the GNN architecture, the learning framework, and the evaluation protocol used throughout the experiments.

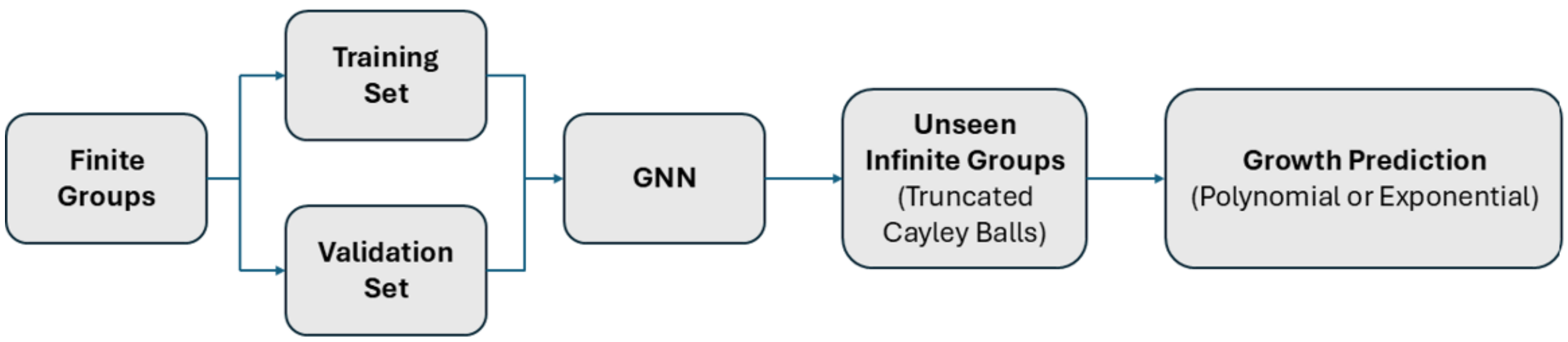

**Figure 1.** Overview of the proposed experimental framework for predicting the growth of finitely generated groups using graph neural networks.

**Datasets.** The experimental design separates the finite training setting from the infinite evaluation setting. Both the training and validation sets consist exclusively of finite groups. The GNN is trained on finite Cayley graphs labeled according to nilpotency, where finite nilpotent groups form the positive class and finite non-nilpotent groups form the negative class. Model selection is performed using only the finite validation set. After the architecture has been selected, the final model is evaluated on truncated Cayley balls ($R = 8$) of unseen infinite groups. In the infinite setting, the learned distinction is interpreted through group growth: virtually nilpotent groups have polynomial growth by Gromov's theorem, whereas the non-virtually nilpotent examples considered here have exponential growth. The following table summarizes the training, validation, and evaluation datasets.

**Table 1.** Groups used in the experimental evaluation.

| Phase | Classification | Group Families | Growth | Purpose |
|---|---|---|---|---|
| Training | Nilpotent | $C_{17}^n \;\; n = 1, 2, 3$<br>$C_{18}^4, C_{20}^4$<br>$H(F_5)$ | | Finite groups used to optimize the model parameters; nilpotent and non-nilpotent groups form the positive and negative training classes, respectively. |
| | Non nilpotent | $SL(2,p) \;\; p = 101, 103, 107$<br>$A_5, S_6$ | | |
| Validation | Nilpotent | $C_{20}, C_{19}^2, C_{21}^3, C_{16}^4$<br>$C_9 \times C_{10} \times C_{10} \times C_9$<br>$C_8 \times C_{12} \times C_{12} \times C_8$<br>$H(F_p) \;\; p = 7,11,13,17,19,23$ | | Finite groups used for model selection and hyperparameter tuning; nilpotent and non-nilpotent groups form the positive and negative validation classes, respectively. |
| | Non nilpotent | $SL(2,p) \;\; p = 109, 127, 131, 149$<br>$S_7, A_6$ | | |

| Testing | Virtually nilpotent | $\mathbb{Z}, \mathbb{Z}^2, \mathbb{Z}^3, \mathbb{Z}^4, \mathbb{Z}^5$<br>$H_3(\mathbb{Z}), D_\infty$<br>$\mathbb{Z}^2 \times C_3, \mathbb{Z}^2 \times C_4$ | Polynomial | Unseen infinite groups used exclusively for final evaluation. Virtually nilpotent and non-virtually nilpotent groups correspond to polynomial and exponential growth, respectively. |
|---|---|---|---|---|
| | Non virtually nilpotent | $F_2, F_3, F_4$<br>$F_2 \times \mathbb{Z}, F_2 \times \mathbb{Z}^2$<br>$F_2 \times C_n \;\; n = 2,3,4,5$<br>$F_2 \times S_3, F_3 \times C_3$ | Exponential | |

The notation used in the table follows standard group-theoretic conventions. The symbols $C_n$ and $C_n^k$ denote finite cyclic groups and finite direct products of $k$ cyclic groups, respectively. The notation $H(F_p)$ denotes the finite Heisenberg group over the finite field with $p$ elements. The groups $SL(2,p)$, $S_n$ and $A_n$ denote the special linear, symmetric, and alternating groups, respectively. The testing dataset consists of the free abelian groups $\mathbb{Z}^k$, the discrete Heisenberg group $H_3(\mathbb{Z})$, the infinite dihedral group $D_\infty$, free groups $F_n$, and direct products involving free groups, free abelian groups, finite cyclic groups, and the symmetric group $S_3$.

**Graph representation**. Each Cayley graph is represented in a format suitable for processing by a GNN. Vertices correspond to group elements, while edges are induced by multiplication by the generators of the group. Each vertex is assigned a two-dimensional feature vector consisting of its degree and its normalized shortest-path distance from the identity element. The graph connectivity is encoded using the edge-index representation of PyTorch Geometric, and the resulting graph serves as the input to the GNN.

**Architecture.** We use a Graph Convolutional Network (GCN) consisting of two graph convolution layers. The first GCNConv layer maps the two-dimensional node features to a hidden representation, and the second GCNConv layer maps this representation to another hidden representation of the same dimension. A ReLU activation is applied after each graph convolution layer. The resulting node embeddings are aggregated using global mean pooling to obtain a graph-level embedding. This embedding is passed through a fully connected hidden layer with ReLU activation, followed by a final linear output layer with two output units corresponding to the two growth classes.

**Learning framework**. The proposed framework is formulated as a binary graph classification problem. The model is trained using supervised learning on the finite training dataset and selected based on its performance on the finite validation dataset. The final model is then applied directly, without retraining, to truncated Cayley balls of radius $R = 8$ of previously unseen infinite finitely generated groups. The experiments therefore evaluate the ability of the GNN to generalize from finite to infinite groups using only finite graph representations.

**Evaluation metric.** Model selection is performed on the finite validation dataset using classification accuracy, which assigns equal importance to both classes regardless of class size. After selecting the best-performing architecture, the final model is evaluated on the previously unseen infinite groups. The evaluation reports the predicted class for each test group together with the overall classification accuracy on the infinite test set.

## 4. Results

This section presents the experimental results obtained using the proposed GNN framework. We first evaluate the candidate models on the finite validation dataset to select the best-performing architecture. The selected model is then evaluated on a collection of previously unseen infinite groups, including both virtually nilpotent and non-virtually nilpotent examples, to assess its ability to generalize beyond the finite groups used during training. The results are analyzed in terms of classification accuracy and the model's ability to distinguish between groups of polynomial and exponential growth.

**Table 2.** Validation accuracy of the evaluated GNN architectures.

| Model | Architecture | Validation accuracy |
|---|---|---|
| 1 | $(2,8,2)$ | $62.5\%$ |
| 2 | $(2,16,2)$ | $83\%$ |
| 3 | $(2,32,2)$ | $72\%$ |
| 4 | $(2,64,2)$ | $83\%$ |

As shown in the table, the architectures with $16$ and $64$ hidden units achieved the highest validation accuracy of 83%. Since both architectures performed equally well on the validation dataset, the smaller architecture, $(2, 16, 2)$, was selected for the subsequent experiments due to its lower complexity.

After selecting the architecture $(2,16,2)$ based on its validation performance, we evaluated the resulting model on the previously unseen infinite groups. This evaluation assesses whether the representations learned exclusively from finite groups generalize to the infinite setting without any additional training. The following table presents the classification results for each test group, including the true and predicted labels together with the corresponding predicted probabilities. The overall classification accuracy on the infinite test set is discussed below.

**Table 3.** Classification results on the unseen infinite groups. Label 1 denotes polynomial growth, whereas label 0 denotes exponential growth. Predicted probabilities are reported in the order (0, 1).

| Group | True | Predicted | P(0) | P(1) |
|---|---|---|---|---|
| $\mathbb{Z}$ | 1 | 1 | 0.23 | 0.77 |
| $\mathbb{Z}^2$ | 1 | 1 | 0.24 | 0.76 |
| $\mathbb{Z}^3$ | 1 | 1 | 0.28 | 0.72 |
| $\mathbb{Z}^4$ | 1 | 0 | 0.53 | 0.47 |
| $\mathbb{Z}^5$ | 1 | 0 | 0.53 | 0.47 |
| $D_\infty$ | 1 | 1 | 0.23 | 0.77 |
| $H_3(\mathbb{Z})$ | 1 | 1 | 0.37 | 0.63 |
| $\mathbb{Z}^2 \times C_3$ | 1 | 1 | 0.25 | 0.75 |
| $\mathbb{Z}^2 \times C_4$ | 1 | 1 | 0.26 | 0.74 |
| $F_2$ | 0 | 0 | 0.98 | 0.02 |
| $F_2 \times \mathbb{Z}$ | 0 | 0 | 1.00 | 0.00 |
| $F_2 \times \mathbb{Z}^2$ | 0 | 0 | 1.00 | 0.00 |
| $F_3$ | 0 | 0 | 1.00 | 0.00 |
| $F_2 \times C_2$ | 0 | 0 | 0.99 | 0.01 |
| $F_2 \times C_3$ | 0 | 0 | 1.00 | 0.00 |
| $F_2 \times C_4$ | 0 | 0 | 1.00 | 0.00 |
| $F_2 \times C_5$ | 0 | 0 | 1.00 | 0.00 |
| $F_2 \times S_3$ | 0 | 0 | 1.00 | 0.00 |
| $F_3 \times C_3$ | 0 | 0 | 1.00 | 0.00 |
| $F_4$ | 0 | 0 | 0.98 | 0.02 |

As shown in the table, the selected GNN model correctly classified 18 of the 20 unseen infinite groups, achieving a test accuracy of 90%. All groups of exponential growth were classified correctly, including the free groups $F_2$, $F_3$, and $F_4$, the direct products $F_2 \times \mathbb{Z}$ and $F_2 \times \mathbb{Z}^2$, the direct products $F_2 \times C_n$ $(n = 2,3,4,5)$, $F_2 \times S_3$, and $F_3 \times C_3$. Among the polynomial-growth groups, the model correctly identified the free abelian groups $\mathbb{Z}$, $\mathbb{Z}^2$ and $\mathbb{Z}^3$, the infinite dihedral group $\mathrm{D}_\infty$, the discrete Heisenberg group $\mathrm{H}_3(\mathbb{Z})$, and the direct products $\mathbb{Z}^2 \times \mathrm{C}_3$ and $\mathbb{Z}^2 \times \mathrm{C}_4$. The only misclassified examples were the higher-rank free abelian groups $\mathbb{Z}^4$ and $\mathbb{Z}^5$, both of which were assigned probabilities close to the decision boundary. These results indicate that the proposed framework successfully generalizes to a broad range of unseen infinite groups while exhibiting limitations for higher-rank free abelian groups.

## 5. Discussion and Conclusion

This work investigated whether GNNs trained exclusively on finite groups can generalize to previously unseen infinite finitely generated groups. The experimental results suggest that such generalization is possible across a diverse collection of infinite groups. Using a validation dataset of finite groups, the architecture $(2,16,2)$ was selected and subsequently evaluated on the infinite test set, achieving a classification accuracy of 90%.

The model correctly classified all non-virtually nilpotent groups of exponential growth and the majority of virtually nilpotent groups of polynomial growth. The only misclassified examples were the higher-rank free abelian groups $\mathbb{Z}^4$ and $\mathbb{Z}^5$. Their predicted probabilities were close to the decision boundary, suggesting that these groups are more challenging to distinguish than the remaining test examples. Overall, the results indicate that the learned graph representations capture structural information related to group growth that generalizes beyond the finite groups used during training.

The present study has several limitations. Although the training and validation datasets include diverse families of finite groups, they remain relatively small and cannot represent the full diversity of finitely generated groups. In addition, the experiments considered a single graph representation and a single GNN architecture while focusing on a binary classification task distinguishing polynomial from exponential growth. Consequently, the reported results should be regarded as an initial demonstration of the proposed framework rather than a comprehensive evaluation of GNNs for learning the growth properties of finitely generated groups.

The proposed framework also depends on the graph representation used to describe each group. In particular, the local structure of a Cayley graph depends on the choice of generating set, even though the growth type of a finitely generated group is independent of that choice. Furthermore, both finite and infinite groups are represented by truncated Cayley balls of fixed radius $R = 8$, so the GNN has access only to local geometric information within this neighborhood of the identity element. Consequently, different choices of generating sets or truncation radii may lead to different graph representations and, therefore, different predictive performance. Investigating the robustness of the proposed framework with respect to these choices remains an important direction for future work.

Future work will investigate larger and more diverse collections of finite and infinite groups, alternative graph representations and GNN architectures, different generating sets and truncation radii, and additional algebraic properties. It would also be of interest to better understand the theoretical relationship between the geometry of Cayley graphs and the ability of graph neural networks to infer global algebraic properties from finite local neighborhoods.